\documentclass[12pt]{article}

\usepackage{amsmath,amsthm,amssymb}
\usepackage{tikz}
\usepackage[colorlinks,
            linkcolor=red,
            anchorcolor=blue,
            citecolor=blue
            ]{hyperref}
\usepackage{pdfsync}
\usepackage[numbers,sort&compress]{natbib}
\allowdisplaybreaks

\theoremstyle{definition}

\theoremstyle{plain}
\newtheorem{theorem}{Theorem}[section]

\newtheorem{remark}{Remark}[section]
\newtheorem{lemma}{Lemma}[section]

\newtheorem{proposition}{Proposition}[section]

\numberwithin{equation}{section}

\begin{document}
	\title{\bf Harmonic functions and the topology of 3-manifolds  under nonnegative scalar curvature }

 \author{Jiannan Chen$^{1}$\footnote{Email address: chenjiannan@email.ncu.edu.cn}, Haiping Fu$^{1}$\footnote{Email address: mathfu@126.com\,\,\    
 	}\\
 	{\small 1.Department of Mathematics,  Nanchang University,} {\small Nanchang 330031, Jiangxi, China.} 
    \\
}

\date{}

\maketitle
\begin{abstract}
In this paper, we establish and recover some results about noncompact 3-manifolds with nonnegative scalar curvature and some topological assumptions. Firstly, we apply the work in Oronzio \cite{Oro1,Oro2} to obtain a nonexistence result of the harmonic function. Secondly, we use the argument in Yan-Zhu \cite{YZ} to establish a rigidity result, which states that such manifold is isometric to the warped product of a flat tours and a half line. Finally, we recover a topological result in Yan-Zhu \cite{YZ} via a different argument.
\end{abstract}	

	

\section{\bf Introduction and main results}
Harmonic functions play a important role in Riemannian geometry. A significant result is Colding's work \cite{Cold1}, where new monotonicity formulas are derived along the level sets of harmonic functions via gradient estimates of the Green's function in nonparbolic Riemannian manifolds with nonnegative Ricci curvature. Subsequently, using the same monotone quantities, Agostiniani-Fogagnolo-Mazzieri \cite{Fog1} obtained the new Willmore-type geometrical inequalities. Furthermore, in \cite{Ben},  Benatti-Fogagnolo-Mazzieri generalized the monotone quantities in \cite{Cold1} to establish an optimal version of the Minskowski inequality in nonparbolic Riemannian manifolds with nonnegative Ricci curvature and Euclidean volume growth. 

In recent years, the level sets of harmonic functions have been used to investigate 3–manifolds with nonnegative scalar curvature. A pioneering work of Stern \cite{Stern} discovered intriguing analogy between the use of such level sets and the use of stable minimal surfaces instituted by Schoen and Yau \cite{SY79}.  Later, Munteanu-Wang \cite{Wang1} established a monotonicity formula for 3-manifolds with nonnegative scalar curvature and provided a sharp comparison results. Then, in \cite{Wang2}, Munteanu-Wang used the harmonic functions to investigate the bottom spectrum of 3-manifolds with scalar curvature lower bound. Furthermore, Chodosh-Li \cite{Otis} generalized Munteanu-Wang's monotonicity formula to 3-manifolds that do not have nonnegative scalar curvature, using this, they resolved the stable Bernstein problem in $ \mathbb{R}^{4} $. Recently, Colding-Minicozzi \cite{Cold4} established the sharp gradient integral estimates to obtain the sharp area comparsion, under the assumption that the norm of gradient vanishes at infinity.  Shortly after, Oronzio \cite{Oro2} generalized the work of Colding–Minicozzi \cite{Cold4} to 3–manifolds with a compact, connected, minimal boundary and established several sharp comparsion results. In \cite{YZ},  Yan–Zhu answered two open questions posed by Wang and Gromov, respectively. More precisely, they investigated the following two problems:
\\
$ \textbf{Wang's question} $ \cite{WJ1,WJ2}. If a contractible 3-manifold admits a complete Riemannian metric with nonnegative scalar curvature, must it be diffeomorphic to $ \mathbb{R}^{3} $?
\\
$ \textbf{Grmov's question} $ \cite{Gro}. Let $ M_{\gamma} $ be the interior of a compact handlebody of genus $ \gamma $. If $ M_{\gamma} $ admits a complete Riemannian metric with
nonnegative scalar curvature, must one have $ \gamma\leq 1 $? 

Under the assumptions that the manifolds admit a lower Ricci curvature bound, and that the minimal positive Green function to vanish at infinity in the nonparabolic case. Yan-Zhu answered both questions affirmatively via a harmonic level-set exhaustions. In the parabolic case, they developed a parabolic version of gradient integral estimates parallel to \cite{Cold4} to obtain the exhaustions whose relevant boundary components are spheres or tori. In the nonparabolic case, they applied the arguments in Colding–Minicozzi \cite{Cold4} to obtain the exhaustions whose relevant boundary components are spheres. Recently, using the Dirichlet Green
function method and the metric replacement technique, You-Zhang \cite{You Z} answered both questions affirmatively without any auxiliary hypothesis. For more relevant results, refer to \cite{WJ1,WJ2,WJ3,Otis1}.

In this paper, we first establish two results about noncompact 3-manifolds with minimal boundary under nonnegative scalar curvature and some  topological assumptions, then we recover a topological result in Yan-Zhu \cite{YZ} via a different argument. Firstly, we apply the work in Oronzio \cite{Oro1,Oro2} to obtain a nonexistence result of the harmonic function. Secondly, we use the argument in Yan-Zhu \cite{YZ} to establish a rigidity result, which states that such manifold is isometric to the warped product of a flat tours and a half line. Finally, we recover a topological result in Yan-Zhu \cite{YZ} via a different argument, which is inspired by the work in \cite{Ag3}, \cite{Oro2} and \cite{You Z}.

We state the main results. First we prove the following nonexistence result of the harmonic function on noncompact 3-manifolds with minimal boundary under nonnegative scalar curvature and some topological assumptions.
\begin{theorem}\label{Th1}
Let $ (M,g) $ be a 3-dimensional, complete, noncompact and one-ended Riemannian manifold with a compact, connected minimal boundary $ \Sigma $ and with nonnegative scalar curvature. Suppose that $ H_{2}(M) $ contains no spherical classes. If the Ricci curvature of $ M $ is bounded from below by a constant, then there exists no a smooth solution $ u $ satisfies  
\begin{equation}\label{u1}
\begin{cases}
\Delta u=0 &\text{in } M,\\
u=0, &\text{on } \Sigma,\\
u\to 1, &\text{at } \infty.
\end{cases}
\end{equation}
\end{theorem}

\begin{remark}
	The topological assumption that $ H_{2}(M) $ contains no spherical classes is motivated by the work of Munteanu-Wang \cite{Wang2}.
\end{remark}

\begin{remark}
In \cite{Wang3},  Munteanu-Wang obtained a nonexistence result about positive Green's function under the assumption that the scalar curvature is bounded from below by a positive constant. Therefore, Theorem~\ref{Th1} can be viewed as an analogue of Munteanu-Wang's result in the setting of the manifolds with minimal boundary and with nonnegative scalar curvature. 
\end{remark}
In \cite{YZ}, Yan-Zhu developed a parabolic version of gradient integral estimates to obtain some topological applications. Inspired by their work, we establish the following rigidity result.
\begin{theorem}\label{Th2}
	Let $ (M,g) $ be a 3-dimensional, complete, noncompact and one-ended Riemannian manifold with a compact, connected minimal boundary $ \Sigma $ and with nonnegative scalar curvature. Suppose that $ H_{2}(M) $ contains no spherical classes. If the Ricci curvature of $ M $ is bounded from below by a constant and there exists a smooth solution $ u $ satisfies  
	\begin{equation}\label{u2}
	\begin{cases}
	\Delta u=0 &\text{in } M,\\
	u=0, &\text{on } \Sigma,\\
	u\to \infty, &\text{at } \infty.
	\end{cases}
	\end{equation}
	Then $ (M, g) $ is isometric to $ \mathbb{T}^{2}\times [0,\infty) $, where $ \mathbb{T}^{2} $ is a tours equipped with a flat metric.
\end{theorem}
In \cite{YZ}, Yan-Zhu obtained the following topological result via the arguments in Colding–Minicozzi \cite{Cold4}. Inspired by the work in \cite{Ag3}, \cite{Oro2} and \cite{You Z}, we give a different proof to recover this result. 
\begin{theorem}\cite{YZ}\label{YZ2}
Let $ (M,g) $ be a complete Riemannian 3-manifold with one end, finite first Betti number and $ R_{g}\geq 0 $. Suppose that $ (M,g) $ is nonparbolic and the Ricci curvature of $ M $ is bounded from below by a constant. If there exists a minimal Green function which vanishes at infinity, then $ M $ admits an exhaustion $ \Omega_{i} $ such that for all $ i $, $ \partial \Omega_{i}\cong \mathbb{S}^{2} $.
\end{theorem}

This paper is organized as follows. In Section 2, we treat Theorem~\ref{Th1}. We first apply the work in Oronzio \cite{Oro1,Oro2} to derive differential inequalities, then establish a quadratic growth of the gradient integral estimates and complete the contradiction argument. In Section 3, we proof Theorem~\ref{Th2}.
 We first use the argument in Yan-Zhu \cite{YZ} to derive a monotonicity formula, then use the minimality of the boundary to obtain the rigidity result. In Section 4, we recover a topological result in Yan-Zhu \cite{YZ}.

\section{\bf Proof of Theorem 1.1}
In this section, we aim to prove Theorem~\ref{Th1}. We first use the work in Oronzio \cite{Oro1,Oro2} to derive two differential inequalities, then establish a quadratic growth of the gradient integral estimates and complete the contradiction argument. To proceed, we first introduce several preliminaries.

Let $ (M,g) $ be a 3-dimensional, complete, noncompact and one-ended Riemannian manifold with a compact and connected boundary $ \Sigma $ and with nonnegative scalar curvature.  Let $ u $ be the solution to (\ref{u2}). We denote the level set of $ u $ 
$$ 
\Sigma_{t}=\{x\in M| u(x)=f(t)\}, 
$$ 
where 
$$ f(t)=1-\frac{\mathfrak{c}_{\Sigma}}{t+\frac{\mathfrak{c}_{\Sigma}}{2}},\,\,\,\,t\in [\frac{\mathfrak{c}_{\Sigma}}{2}, +\infty).
$$
Set
$$
\mathcal{T}=\Big\{t\in [\frac{\mathfrak{c}_{\Sigma}}{2}, +\infty)\,\,\Big|\,\, f(t) \text{ is a regular value of } u \Big\}
$$
Let $ \mathfrak{c}_{\Sigma}>0 $ be the boundary capacity of $ \Sigma $ in $ (M,g) $, that is 
$$
\mathfrak{c}_{\Sigma}=\frac{1}{4\pi}\int_{M}|\nabla u|^{2}=\frac{1}{4\pi}\int_{\Sigma}|\nabla u|.
$$
Along the level set $ \Sigma_{t} $, we define 
$$
\mathcal{A}(t)=\frac{t^{2}}{\mathfrak{c}_{\Sigma}^{2}}\Big(1+\frac{\mathfrak{c}_{\Sigma}}{2t}\Big)^{4}\int_{\Sigma_{t}}|\nabla u|^{2},
$$
$$
\mathcal{G}(t)=-\frac{\pi\mathfrak{c}_{\Sigma}^{2}}{t}+\frac{t}{4}\Big(1+\frac{\mathfrak{c}_{\Sigma}}{2t} \Big)^{4}\int_{\Sigma_{t}}|\nabla u|^{2},
$$
$$
\mathcal{F}(t)=4\pi t+\frac{t^{3}}{\mathfrak{c}_{\Sigma}^{2}}\Big(1+\frac{\mathfrak{c}_{\Sigma}}{2t}\Big)^{3}\Big(1-\frac{3\mathfrak{c}_{\Sigma}}{2t}\Big)\int_{\Sigma_{t}}|\nabla u|^{2}-\frac{t^{2}}{\mathfrak{c}_{\Sigma}}\Big(1+\frac{\mathfrak{c}_{\Sigma}}{2t}\Big)^{2}\int_{\Sigma_{t}}|\nabla u|H.
$$
Set
$$
A_{0}=\mathcal{F}(\frac{\mathfrak{c}_{\Sigma}}{2})=-\mathcal{G}(\frac{\mathfrak{c}_{\Sigma}}{2})=2\mathfrak{c}_{\Sigma}\Big(\pi-\int_{\Sigma}|\nabla u|^{2}\Big),
$$
\begin{equation}\label{A1 A0}
A_{1}=2\pi\mathfrak{c}_{\Sigma}- A_{0}=2\mathfrak{c}_{\Sigma}\int_{\Sigma}|\nabla u|^{2}.
\end{equation}
By Hopf's lemma, $ 0 $ is a regular value, hence 
\begin{equation}\label{A1 geq 0}
A_{1}>0.
\end{equation}
We also denote
$$
\widetilde{\mathcal{A}}(t)=\mathcal{A}(t)-\frac{1}{2t}A_{1}.
$$
As in \cite{Oro2}, we introduce the $ (0,2) $-symmetric tensors of $ M $:
$$
T=\nabla^{2}u+\frac{6u}{1-u^{2}}du\otimes du,\,\,\,B=T-\frac{tr(T)}{3}g,
$$
We define
$$
B_{1}(t)=\int_{\Sigma_{t}}\frac{|B|^{2}}{|\nabla u|^{2}},\,\,\,R_{1}(t)=\int_{\Sigma_{t}}R_{g}.
$$
Observe that 
$$
\frac{|B|^{2}}{2|\nabla u|^{2}}=\frac{|\nabla^{\Sigma_{t}}|\nabla u| |^{2}}{|\nabla u|^{2}}+\frac{1}{2}|\overset{\circ}{\mathbb{II}}|^{2}+\frac{3}{4}\Big(H-\frac{4u}{1-u^{2}}|\nabla u|\Big)^{2}.
$$
In \cite{Oro1}, Oronzio considered $ \mathcal{F}(t) $ for harmonic functions. Hence, we have the following result.
\begin{proposition}\cite{Oro1}\label{F mono}
	Let $ (M,g) $ be a complete, orientable, one-ended 3-manifold with compact, connnected boundary $ \Sigma $. Suppose that $ H_{2}(M,\Sigma)=0 $. Let $ u $ be the solution of (\ref{u1}). If $ R_{g}\geq 0 $, then the following quantity 
	\begin{equation}\label{F}
\mathcal{F}(t)=4\pi t+\frac{t^{3}}{\mathfrak{c}_{\Sigma}^{2}}\Big(1+\frac{\mathfrak{c}_{\Sigma}}{2t}\Big)^{3}\Big(1-\frac{3\mathfrak{c}_{\Sigma}}{2t}\Big)\int_{\Sigma_{t}}|\nabla u|^{2}-\frac{t^{2}}{\mathfrak{c}_{\Sigma}}\Big(1+\frac{\mathfrak{c}_{\Sigma}}{2t}\Big)^{2}\int_{\Sigma_{t}}|\nabla u|H
	\end{equation}
	is monotone nondecreasing in $ t $. Futhermore, at a regular value $ t $, we have 
	\begin{equation}\label{Fk}
	\begin{aligned}
	\mathcal{F}^{\prime}(t)&=4\pi-2\pi\chi(\Sigma_{t})+\frac{1}{2}\int_{\Sigma_{t}}\Big(\frac{|B|^{2}}{|\nabla u|^{2}}+R_{g} \Big).
	\end{aligned}
	\end{equation}
	Moreover, if $ \mathcal{F}(t)\equiv0 $ for all  $ t\in [\frac{\mathfrak{c}_{\Sigma}}{2}, +\infty) $, then $ (M,g) $ is isometric to a spatial Schwarzschild manifolds outside a rotationally symmetric sphere.
\end{proposition}
Using the monotone quantity $ \mathcal{F}(t) $, one has the following differential inequalities.
\begin{proposition}\label{diff ineq}
	Let $ (M,g) $ be a complete, orientable, one-ended 3-manifold with compact, connnected boundary $ \Sigma $ and with nonnegative scalar curvature. Suppose that  $ H_{2}(M) $ contains no spherical classes and $ H_{2}(M,\Sigma)=0 $. Let $ u $ be the solution of (\ref{u1}). Then we have the following differential inequalities.
	\\
	$ (i) $ For all $ t\in [\frac{\mathfrak{c}_{\Sigma}}{2}, +\infty) $,
	$$
	t\mathcal{A}^{\prime}(t)\geq \mathcal{A}(t)-\frac{A_{1}}{t}+\frac{1}{2t}\int_{\frac{\mathfrak{c}_{\Sigma}}{2}}^{t}(B_{1}(s)+R_{1}(s))ds
	$$
	holds everywhere.
	\\
	$ (ii) $ For all $ t\in [\frac{\mathfrak{c}_{\Sigma}}{2}, +\infty) $, 
	$$
	t\widetilde{\mathcal{A}}^{\prime}(t)\geq \widetilde{\mathcal{A}}(t)+\frac{1}{2t}\int_{\frac{\mathfrak{c}_{\Sigma}}{2}}^{t}(B_{1}(s)+R_{1}(s))ds
	$$
	holds everywhere.
\end{proposition}
\begin{proof}
At a regular value $ t\in \mathcal{T} $, we compute
$$
\mathcal{A}^{\prime}(t)=\frac{2t}{\mathfrak{c}_{\Sigma}^{2}}\Big(1+\frac{\mathfrak{c}_{\Sigma}}{2t}\Big)^{3}\Big(1-\frac{\mathfrak{c}_{\Sigma}}{2t}\Big)\int_{\Sigma}|\nabla u|^{2}-\frac{1}{\mathfrak{c}_{\Sigma}}\Big(1+\frac{\mathfrak{c}_{\Sigma}}{2t}\Big)^{2}\int_{\Sigma_{t}}|\nabla u|H.
$$
Then we get 
\begin{equation}\label{t Ak}
\begin{aligned}
t\mathcal{A}^{\prime}(t)-\mathcal{A}(t)=&\frac{t^{2}}{\mathfrak{c}_{\Sigma}^{2}}\Big(1+\frac{\mathfrak{c}_{\Sigma}}{2t}\Big)^{3}\Big(1-\frac{3\mathfrak{c}_{\Sigma}}{2t}\Big)\int_{\Sigma}|\nabla u|^{2}\\
&-\frac{t}{\mathfrak{c}_{\Sigma}}\Big(1+\frac{\mathfrak{c}_{\Sigma}}{2t}\mathfrak{c}_{\Sigma}\Big)^{2}\int_{\Sigma_{t}}|\nabla u|H\\
=&\frac{1}{t}\mathcal{F}(t)-4\pi\\
=&\frac{1}{t}\Big(\mathcal{F}(\frac{\mathfrak{c}_{\Sigma}}{2})+\int_{\frac{\mathfrak{c}_{\Sigma}}{2}}^{t} \mathcal{F}^{\prime}(s)ds\Big)-4\pi.
\end{aligned}
\end{equation}
On the other hand, since $ H_{2}(M) $ contains no spherical classes, all connected components of $ \Sigma_{t} $ are not spheres, by the Guass-Bonnet theorem and (\ref{Fk}) , we have
\begin{equation}
\begin{aligned}\label{F ineq}
\mathcal{F}^{\prime}(t)&=4\pi-2\pi\chi(\Sigma_{t})+\frac{1}{2}\int_{\Sigma_{t}}\Big(\frac{|B|^{2}}{|\nabla u|^{2}}+R_{g} \Big)\\
&\geq 4\pi+\frac{1}{2}(B_{1}(t)+R_{1}(t)).
\end{aligned}
\end{equation}
Thus we find 
\begin{equation}\label{diff Ak}
\begin{aligned}
t\mathcal{A}^{\prime}(t)-\mathcal{A}(t)&\geq \frac{1}{t}\mathcal{F}(\frac{\mathfrak{c}_{\Sigma}}{2})+\frac{4\pi(t-\frac{\mathfrak{c}_{\Sigma}}{2})}{t}-4\pi+\frac{1}{2t}\int_{\frac{\mathfrak{c}_{\Sigma}}{2}}^{t}(B_{1}(s)+R_{1}(s))ds\\
&=\frac{1}{t}\Big(A_{0}-2\pi \mathfrak{c}_{\Sigma}\Big)+\frac{1}{2t}\int_{\frac{\mathfrak{c}_{\Sigma}}{2}}^{t}(B_{1}(s)+R_{1}(s))ds\\
&=-\frac{A_{1}}{t}+\frac{1}{2t}\int_{\frac{\mathfrak{c}_{\Sigma}}{2}}^{t}(B_{1}(s)+R_{1}(s))ds.
\end{aligned}
\end{equation}
By Sard's theorem and continuity of $ \mathcal{A}(t) $, (\ref{diff Ak}) holds for every $ t\in [\frac{\mathfrak{c}_{\Sigma}}{2}, +\infty) $. This proves $ (i) $.

To prove $ (ii) $, a simple calculation shows that
$$
\begin{aligned}
t\widetilde{\mathcal{A}}^{\prime}(t)-\widetilde{\mathcal{A}}(t)&=t\mathcal{A}^{\prime}(t)-\mathcal{A}(t)+\frac{1}{t}A_{1}\geq \frac{1}{2t}\int_{\frac{\mathfrak{c}_{\Sigma}}{2}}^{t}(B_{1}(s)+R_{1}(s))ds.
\end{aligned}
$$
This proves $ (ii) $.
\end{proof}
Using Proposition~\ref{diff ineq}, we now derive the quadratic growth of $ \mathcal{A}(t) $. 
\begin{proposition}\label{quadratic growth}
	Under the assumptions in Proposition~\ref{diff ineq}, if $ \mathcal{A}(t_{0})>\frac{A_{1}}{t_{0}} $ for some $ t_{0}\in [\frac{\mathfrak{c}_{\Sigma}}{2},+\infty ) $, then there exists a positive constant $ C>0 $ such that, for every $ t\geq t_{0} $,
	$$
	\mathcal{A}(t)\geq Ct^{2}.
	$$
\end{proposition}
\begin{proof}
First we choose $ \delta>0 $ such that 
\begin{equation}\label{wAk t0}
\widetilde{\mathcal{A}}(t_{0})>\delta.
\end{equation}
In fact, since $ \mathcal{A}(t_{0})>\frac{A_{1}}{t_{0}}>0 $ by (\ref{A1 geq 0}), we have 
$$ \widetilde{\mathcal{A}}(t_{0})=\mathcal{A}(t_{0})-\frac{1}{2t_{0}}A_{1}>\frac{A_{1}}{2t_{0}}>0.
$$
Thus, we can choose $ \delta>0 $ small enough such that (\ref{wAk t0}) holds.
We claim that there exists $ \delta>0 $, such that for every $ t\geq t_{0} $,
\begin{equation}\label{key estimate}
\mathcal{A}^{\prime}(t)\geq 0,\,\,\,\,\mathcal{A}(t)\geq \delta\frac{t}{t_{0}}.
\end{equation}
Indeed, let $ \delta>0 $ be chosen as in (\ref{wAk t0}). By Proposition~\ref{diff ineq} $ (i) $, we have 
\begin{equation}\label{diff ineq case2}
t\mathcal{A}^{\prime}(t)\geq \mathcal{A}(t)-\frac{A_{1}}{t}.
\end{equation}
Since $ \mathcal{A}(t_{0})>4\pi-\frac{A_{1}}{t_{0}} $, one has $ \mathcal{A}^{\prime}(t_{0})>0 $. This shows that 
$$
\mathcal{A}(t)>\frac{A_{1}}{t}
$$
for all $ t\geq t_{0} $, which implies that $ \mathcal{A}^{\prime}(t)\geq0 $ for all $ t\geq t_{0} $. Indeed, we consider the function $ f(t)=\frac{1}{t}(\mathcal{A}(t)-\frac{A_{1}}{t}) $, then (\ref{diff ineq case2}) gives 
$$
f^{\prime}(t)=\frac{1}{t^{2}}\Big(t\mathcal{A}^{\prime}(t)-\mathcal{A}_{k}(t)+\frac{2A_{1}}{t} \Big)\geq \frac{1}{t^{2}}\Big(t\mathcal{A}^{\prime}(t)-\mathcal{A}(t)+\frac{A_{1}}{t} \Big)\geq 0.
$$
due to $ A_{1}>0 $. Hence $ f(t)\geq f(t_{0})>0 $, that is $ \mathcal{A}(t)>\frac{A_{1}}{t} $. On the other hand, by Proposition~\ref{diff ineq} $ (ii) $, we have 
\begin{equation}\label{wide A}
t\widetilde{\mathcal{A}}^{\prime}(t)\geq \widetilde{\mathcal{A}}(t).
\end{equation}
Since $ \widetilde{\mathcal{A}}(t_{0})>0 $, we deduce that $ \widetilde{\mathcal{A}}^{\prime}(t_{0})>0 $, and so 
$$
\widetilde{\mathcal{A}}(t)>0,\,\,\,\widetilde{\mathcal{A}}^{\prime}(t)\geq 0,
$$
for all $ t\geq t_{0} $. Indeed, we consider the function $ h(t)=\frac{1}{t}\widetilde{\mathcal{A}}(t) $, then (\ref{wide A}) gives 
$$
h^{\prime}(t)=\frac{1}{t^{2}}\Big(t\widetilde{\mathcal{A}}^{\prime}(t)-\widetilde{\mathcal{A}}(t)\Big)\geq 0.
$$
Hence $ h(t)\geq h(t_{0})>0 $, that is $ \widetilde{\mathcal{A}}(t)>0 $, and so $ \widetilde{\mathcal{A}}^{\prime}(t)\geq0 $. Now (\ref{wide A}) can be rewritten as 
$$
\frac{\widetilde{\mathcal{A}}^{\prime}(t)}{\widetilde{\mathcal{A}}(t)}\geq \frac{1}{t}
$$
for all $ t\geq t_{0} $. Integrating from $ t_{0} $ to $ t $ gives 
$$
\widetilde{\mathcal{A}}(t)\geq \widetilde{\mathcal{A}}(t_{0})\frac{t}{t_{0}}>\delta\frac{t}{t_{0}}
$$
for all $ t\geq t_{0} $. Since $ A_{1}>0 $, we have
$$ \widetilde{\mathcal{A}}(t)=\mathcal{A}(t)-\frac{1}{2t}A_{1}<\mathcal{A}_{k}(t).
$$ 
Therefore, we conclude that 
$$
\mathcal{A}(t)\geq \delta\frac{t}{t_{0}}
$$
for all $ t\geq t_{0} $. Hence (\ref{key estimate}) holds. As in \cite{Cold4,Oro2}, on the interval $ [\frac{\mathfrak{c}_{\Sigma}}{2}, +\infty) $, we define the function 
$$
a(t)=t\frac{\mathcal{A}^{\prime}(t)}{\mathcal{A}(t)}.
$$
From (\ref{key estimate}), we see that $ a(t)\geq0 $ for all $ t\geq t_{0} $. By (\ref{t Ak}), for every $ t\in \mathcal{T} $, one has
$$
t\mathcal{A}^{\prime}(t)-\mathcal{A}(t)=\frac{1}{t}\mathcal{F}(t)-4\pi.
$$ 
This leads to 
$$
\begin{aligned}
\mathcal{A}^{\prime}(t)&=\frac{\mathcal{F}(t)}{t^{2}}+\frac{1}{t}(\mathcal{A}(t)-4\pi),\\
t\mathcal{A}^{\prime\prime}(t)&=[t\mathcal{A}^{\prime}(t)-\mathcal{A}(t)]^{\prime}=-\frac{\mathcal{F}(t)}{t^{2}}+\frac{\mathcal{F}^{\prime}(t)}{t}.
\end{aligned}
$$
Hence one has 
$$
a^{\prime}(t)=\frac{1}{t}\Big(1-\frac{4\pi}{\mathcal{A}(t)}+\frac{\mathcal{F}^{\prime}(t)}{\mathcal{A}(t)}-a^{2}(t) \Big).
$$
For $ \nu=\frac{\nabla u}{|\nabla u|} $, observe that 
$$
\begin{aligned}
B(\nu,\nu)&=-H|\nabla u|+\frac{4u}{1-u^{2}}|\nabla u|^{2}\\
&=-H|\nabla u|+\frac{2}{\mathfrak{c}_{\Sigma}}\Big(1+\frac{\mathfrak{c}_{\Sigma}}{2t}\Big)\Big(1-\frac{\mathfrak{c}_{\Sigma}}{2t} \Big)|\nabla u|^{2}.
\end{aligned}
$$
Then 
$$
t\mathcal{A}^{\prime}(t)=\frac{t}{\mathfrak{c}_{\Sigma}}\Big(1+\frac{\mathfrak{c}_{\Sigma}}{2t}\Big)^{2}\int_{\Sigma_{t}}B(\nu,\nu).
$$
As in \cite{Oro2}, one has 
$$
t^{2}(\mathcal{A}^{\prime}(t))^{2}\leq \frac{2}{3}\mathcal{A}(t)B_{1}(t).
$$
On the other hand, by (\ref{F ineq}) and $ R_{g}\geq 0 $, we have 
$$
\mathcal{F}^{\prime}(t)\geq 4\pi+\frac{1}{2}B_{1}(t).
$$
Thus one has 
$$
\frac{\mathcal{F}^{\prime}(t)}{\mathcal{A}(t)}\geq \frac{4\pi}{\mathcal{A}(t)}+\frac{B_{1}(t)}{2\mathcal{A}(t)}\geq \frac{4\pi}{\mathcal{A}(t)}+\frac{3}{4}a(t)^{2}.
$$
Therefore, 
\begin{equation}\label{a ineq}
a^{\prime}(t)\geq \frac{1}{t}\Big(1-\frac{1}{4}a(t)^{2} \Big)
\end{equation}
for each $ t\in \mathcal{T} $.
Now, we show the existence of a positive constant $ C=C(t_{0}) $ such that 
$$
a(t)\geq 2-\frac{C}{\sqrt{t}}
$$
for all $ t\geq t_{0} $. If $ a(t)\geq 2 $ for all $ t\geq t_{0} $, the claim is trivially true. Therefore, we suppose that there exists $ \bar{t}_{0}\geq t_{0} $ such that 
$$
a(\bar{t}_{0})<2.
$$
Without loss of generality, we may assume that $ \bar{t}_{0}>t_{0} $ by the continuity of $ a(t) $ and that $ \bar{t}_{0}\in \mathcal{T} $ by Sard's theorem.

Let $ t_{2}>t_{0} $ be a generic regular value such that $ a(t_{2})<2 $. Denote 
$$
t_{1}=\inf\limits_{t\in [t_{0},t_{2}]}\{t\,\,|\,\,a(s)<2,\,\,\forall s\in (t,t_{2}] \}.
$$ 
Either $ t_{1}=t_{0} $ or $ t_{1}>t_{0} $. We will consider two cases.

Suppose first that $ t_{1}=t_{0} $. In this case, one has 
\begin{equation}\label{a(t) 0 2}
0\leq a(t)\leq 2,\,\,\,\forall t\in [t_{0}, t_{2}].
\end{equation}
Hence, it follows from (\ref{a ineq}) that 
$$
a^{\prime}(t)\geq \frac{1}{t}\Big(1-\frac{1}{2}a(t) \Big),\,\,\,\,\,\forall t\in [t_{0},t_{2}].
$$
This implies that 
\begin{equation}\label{a ineq1}
(\sqrt{t}a(t))^{\prime}\geq \frac{1}{\sqrt{t}}
\end{equation}
for every $ t\in \mathcal{T}\cap [t_{0},t_{2}] $.
Observe that the function $ \sqrt{t}a(t) $ is locally absolute continuous in $ [\frac{\mathfrak{c}_{\Sigma}}{k},+\infty)$. Thus, integrating (\ref{a ineq1}) from $ t_{0} $ to $ t_{2} $ gives 
$$
a(t_{2})\geq 2+\sqrt{\frac{t_{0}}{t_{2}}}a(t_{0})-2\sqrt{\frac{t_{0}}{t_{2}}}\geq 2-2\sqrt{\frac{t_{0}}{t_{2}}}.
$$ 
In the remaining case, we have  $ t_{1}>t_{0} $. Thus one deduces that $ a(t_{1})=2 $ and that $ 0\leq a(t)\leq 2 $, for any $ t\in [t_{1},t_{2}] $. Hence, proceeding as above, we get 
$$
a(t_{2})\geq 2+\sqrt{\frac{t_{1}}{t_{2}}}a(t_{1})-2\sqrt{\frac{t_{1}}{t_{2}}}=2.
$$
Therefore, choosing $ C=C(t_{0})=2\sqrt{t_{0}} $, one has 
\begin{equation}\label{a regul}
a(t)\geq 2-\frac{C}{\sqrt{t}},
\end{equation}
for every $ t\in \mathcal{T}\cap [t_{0},+\infty) $. The continuity of the function $ a(t) $ implies that (\ref{a regul}) holds everywhere in the interval $ [t_{0},+\infty) $. Then, by the definition of $ a(t) $, we have 
$$
(\log(\mathcal{A}(t)))^{\prime}\geq 2t^{-1}-Ct^{-\frac{3}{2}}
$$ 
for all $ t\in [t_{0},+\infty) $. Integrating this inequality from $ t_{0} $ to $ t $, we obtain that 
$$
\mathcal{A}(t)\geq t_{0}^{-2}\mathcal{A}(t_{0})e^{-\frac{2C}{\sqrt{t_{0}}}}t^{2},
$$
which is the desired quadratic growth of $ \mathcal{A}(t) $.
\end{proof}

\textbf{Proof of Theorem 1.1}. Suppose by contradiction that there exists a solution $ u $ satisfies (\ref{u1}). Then we claim that 
\begin{equation}\label{Ak A1}
\mathcal{A}(t)\leq\frac{A_{1}}{t}
\end{equation}
for all $ t\in [\frac{\mathfrak{c}_{\Sigma}}{2}, +\infty) $. In fact, if there exists a $ t_{0}\in [\frac{\mathfrak{c}_{\Sigma}}{2}, +\infty)  $ such that $ \mathcal{A}(t_{0})>\frac{A_{1}}{t_{0}} $, then 
$$ 
\mathcal{A}(t)\geq Ct^{2} 
$$ 
for all $ t\in [t_{0}, +\infty) $, , which can be rewritten as 
$$
\int_{\Sigma_{t}}|\nabla u|^{2}\geq C\Big(1+\frac{\mathfrak{c}_{\Sigma}}{2t}\Big)^{-4}\mathfrak{c}_{\Sigma}^{2}\geq 2^{-4}C\mathfrak{c}_{\Sigma}^{2}.
$$
However, by the Cheng-Yau gradient estimate \cite{Cheng-Yau}, 
$$
\sup_{\Sigma_{t}}|\nabla u|=\sup_{\Sigma_{t}}|\nabla (1-u)|\leq C(1-f(t))\to 0
$$
as $ t\to +\infty $. Therefore, we get 
$$
\int_{\Sigma_{t}}|\nabla u|^{2}\leq \sup_{\Sigma_{t}}|\nabla u|\int_{\Sigma_{t}}|\nabla u|=\mathfrak{c}_{\Sigma}\sup_{\Sigma_{t}}|\nabla u|\to 0.
$$
This is a contradiction. Hence (\ref{Ak A1}) holds, which implies that 
$$
\lim\limits_{t\to +\infty}t\Big(1+\frac{\mathfrak{c}_{\Sigma}}{2t}\Big)^{4}\int_{\Sigma_{t}}|\nabla u|^{2}=0.
$$
Hence, we obtain that 
\begin{equation}\label{G lim}
\lim\limits_{t\to +\infty}\mathcal{G}(t)=\lim\limits_{t\to +\infty}\Big[-\frac{\pi\mathfrak{c}_{\Sigma}^{2}}{t}+\frac{t}{4}\Big(1+\frac{\mathfrak{c}_{\Sigma}}{2t} \Big)^{4}\int_{\Sigma_{t}}|\nabla u|^{2}\Big]=0.
\end{equation}
By (\ref{F ineq}), we have $ \mathcal{F}^{\prime}(t)\geq 4\pi $, so that
$$
\mathcal{F}(t)\geq \mathcal{F}(\frac{\mathfrak{c}_{\Sigma}}{2})+4\pi(t-\frac{\mathfrak{c}_{\Sigma}}{2}).
$$ 
Thus, for every $ t\in \mathcal{T} $, one has
\begin{equation}\label{G diff}
\mathcal{G}^{\prime}(t)=\frac{\mathfrak{c}_{\Sigma}^{2}}{4t^{3}}\mathcal{F}(t)\geq \frac{\mathfrak{c}_{\Sigma}^{2}}{4t^{3}}\Big(\mathcal{F}(\frac{\mathfrak{c}_{\Sigma}}{2})+4\pi(t-\frac{\mathfrak{c}_{\Sigma}}{2})\Big).
\end{equation}
Integrating (\ref{G diff}) from $ \frac{\mathfrak{c}_{\Sigma}}{2} $ to $ t $ gives 
$$
\mathcal{G}(t)-\mathcal{G}(\frac{\mathfrak{c}_{\Sigma}}{2})\geq \Big(\frac{1}{2}-\frac{\mathfrak{c}_{\Sigma}^{2}}{8t^{2}}\Big)\mathcal{F}(\frac{\mathfrak{c}_{\Sigma}}{2})+\pi\mathfrak{c}_{\Sigma}^{2}\Big(\frac{\mathfrak{c}_{\Sigma}}{4t^{2}}-\frac{1}{t}+\frac{1}{\mathfrak{c}_{\Sigma}}\Big).
$$
Letting $ t\to +\infty $ and using (\ref{G lim}), we get 
\begin{equation}\label{A0 ineq}
A_{0}=\mathcal{F}(\frac{\mathfrak{c}_{\Sigma}}{2})=-\mathcal{G}(\frac{\mathfrak{c}_{\Sigma}}{2})\geq \frac{1}{2}A_{0}+\pi\mathfrak{c}_{\Sigma}.
\end{equation}
On the other hand, by (\ref{A1 A0}), we have 
\begin{equation}\label{A1 0}
A_{1}=2\pi\mathfrak{c}_{\Sigma}- A_{0}=2\mathfrak{c}_{\Sigma}\int_{\Sigma}|\nabla u|^{2}. 
\end{equation}
Hence, it follows frm (\ref{A0 ineq}) and (\ref{A1 0}) that   
$$
A_{1}\leq 0.
$$
This contradicts to (\ref{A1 geq 0}). Therefore, such a solution $ u $ cannot exist. This completes the proof of Theorem~\ref{Th1}.
$\hfill\square$

\section{\bf Proof of Theorem 1.2}
In this section, we aim to prove Theorem~\ref{Th2}. The main argument follws Yan-Zhu \cite{YZ}.  We first use the method in Yan-Zhu \cite{YZ} to derive a monotonicity formula, then use the minimality of the boundary to obtain the rigidity result. To proceed, we first introduce some preliminaries.

Under the assumption that $ H_{2}(M) $ contains no spherical classes,  we know that the level set 
$$ 
\Sigma_{r}=\{u=r\} 
$$
is a smooth and closed surface of genus at least one for every regular value $ r\geq 0 $, namely 
$$ 
\text{ genus}(\Sigma_{r})\geq 1.
$$
We define 
$$
\mathcal{A}_{1}(r)=\int_{\Sigma_{r}}|\nabla u|^{2},\,\,\,\mathcal{B}_{1}(r)=\mathcal{A}_{1}^{\prime}(r).
$$
Let $ \mathbf{n}=\frac{\nabla u}{|\nabla u|} $ be the unit normal of $\Sigma_{r}$.

The following two lemmas follows from a slight modification of the proof of Yan-Zhu \cite{YZ}, together with our assumption that $ H_{2}(M) $ contains no spherical classes.
\begin{lemma}\label{pre}\cite{YZ}
We have the following results.
\\
$ (i) $ $\int_{\Sigma_{r}}|\nabla u|=\text{const}$ for $ r\geq 0 $;
\\
$ (ii) $ $ \mathcal{A}_{1}^{\prime}(r)=\int_{\Sigma_{r}}\nabla^{2}u(\mathbf{n},\mathbf{n})=-\int_{\Sigma_{r}}H|\nabla u|=\mathcal{B}_{1}(r) $;
\\
$ (iii) $ If $ r_{1}<r_{2} $ are regular values, then 
$$
\mathcal{B}_{1}(r_{2})-\mathcal{B}_{1}(r_{1})=\int_{r_{1}}^{r_{2}}\int_{\Sigma_{r}}\frac{1}{2}\Big(\frac{|\nabla^{2} u |^{2}}{|\nabla u|^{2}}+R_{g}-R_{\Sigma_{r}} \Big)\geq 0.
$$
\end{lemma}

\begin{lemma}\label{quad A1}\cite{YZ}
Assume that $ r_{0} $ is a regular value such that 
$$
\mathcal{B}_{1}(r_{0})=\mathcal{A}_{1}^{\prime}(r_{0})>0.
$$
Then there exists a constant $ C>0 $ such that 
$$
\mathcal{A}_{1}(r)\geq Cr^{2}.
$$
for all $ r\geq r_{0} $.
\end{lemma}
\textbf{Proof of Theorem 1.2}. First we show that 
\begin{equation}\label{B1 leq0}
\mathcal{B}_{1}(r)\leq 0
\end{equation}
for all $ r\geq 0 $. Indeed, otherwise, there exists $ r_{0}\geq 0 $ such that $ \mathcal{B}_{1}(r_{0})>0 $. Then, by Lemma~\ref{quad A1}, there exists a constant $ C>0 $ such that 
$$
\mathcal{A}_{1}(r)\geq Cr^{2}.
$$
for $ r\geq r_{0} $. On the other hand, by the Cheng-Yau gradient estimate \cite{Cheng-Yau}, for all sufficiently large $ r $,
$$
\sup_{\Sigma_{r}}|\nabla u|\leq Cr.
$$
Combining this with Lemma~\ref{pre} $ (i) $, we get 
$$
\mathcal{A}_{1}(r)=\int_{\Sigma_{r}}|\nabla u|^{2}\leq \sup_{\Sigma_{r}}|\nabla u|\int_{\Sigma_{r}}|\nabla u|\leq Cr.
$$
This is a contradiction. Therefore, (\ref{B1 leq0}) holds. From Lemma~\ref{pre} $ (iii) $ and (\ref{B1 leq0}), we know that $ \mathcal{B}_{1} $ is monotone nondecreasing and bounded above by $ 0 $, so it admits a limit
$$
l:=\lim\limits_{r\to \infty}\mathcal{B}_{1}(r)\leq 0.
$$
If $ l<0 $, then there exists $ r_{0} $ large enough such that, for all $ r\geq r_{0} $, 
\begin{equation}\label{B l}
\mathcal{A}_{1}^{\prime}(r)=\mathcal{B}_{1}(r)\leq \frac{l}{2}<0.
\end{equation}
Integrating (\ref{B l}) from $ 0 $ to $ r>0 $ gives 
$$
0\leq \mathcal{A}_{1}(r)\leq \mathcal{A}_{1}(0)+\frac{l}{2}r\to -\infty
$$
as $ r\to \infty $, which is impossible. Therefore, 
$$
\lim\limits_{r\to \infty}\mathcal{B}_{1}(r)=0.
$$
On the other hand, notice that the boundary $ \Sigma $ is minimal, we have $$ 
\mathcal{B}_{1}(0)=-\int_{\Sigma}H|\nabla u|=0.
$$ 
Since $ \mathcal{B}_{1} $ is monotone nondecreasing , we deduce that 
$$
\mathcal{B}_{1}(r)\equiv0
$$
for all $ r\geq 0 $. Therefore, by Lemma~\ref{pre} $ (iii) $, one has 
$$
\int_{0}^{\infty}\int_{\Sigma_{r}}\frac{1}{2}\Big(\frac{|\nabla^{2} u |^{2}}{|\nabla u|^{2}}+R_{g}-R_{\Sigma_{r}} \Big)\equiv0.
$$
By Gauss-Bonnet theorem, this can be rewritten as 
$$
\int_{0}^{\infty}\Big(\int_{\Sigma_{r}}\Big(\frac{|\nabla^{2} u |^{2}}{|\nabla u|^{2}}+R_{g}\Big)+4\pi(\text{genus}(\Sigma_{r})-1) \Big)dr\equiv0.
$$
Therefore, $ |\nabla^{2} u|\equiv0 $, $ R_{g}\equiv0 $ and $ \text{genus}(\Sigma_{r})\equiv1 $. Hence $ |\nabla u|\equiv \text{const}\neq0 $, which implies that all $ r\geq0 $ are regular values and so all $ \Sigma_{r} $ are connected and diffeomorphic to $ \mathbb{T}^{2} $. Consequently, $ M $ is isometric to $  \mathbb{T}^{2}\times [0,+\infty) $, where $ \mathbb{T}^{2} $ is a tours equipped with a flat metric.
$\hfill\square$
\section{Exhuastion of spheres}
In this section, we recover the following topological result in Yan-Zhu \cite{YZ}. This argument is inspired by the work in \cite{Ag3}, \cite{Oro2} and \cite{You Z}. Let $ G $ be the minimal Green function with pole at $ p\in M $ in Theorem~\ref{YZ2} and let $ u=1-4\pi G $. We state Theorem~\ref{YZ2} as the following equivalent form. 
\begin{proposition}\label{YZ}\cite{YZ}
Let $ (M,g) $ be a complete Riemannian 3-manifold with one end, finite first Betti number and $ R_{g}\geq 0 $. Suppose that $ (M,g) $ is nonparbolic and the Ricci curvature of $ M $ is bounded from below by a constant. If there exists a maximal distributional solution $ u $ satisfies  
\begin{equation}\label{u3}
\begin{cases}
\Delta u=4\pi \delta_{p} &\text{in } M,\\
u\to 1, &\text{at } \infty.
\end{cases}
\end{equation}
Then $ M $ admits an exhuastion $ \Omega_{i} $ such that for all $ i $, $ \partial \Omega_{i}\cong \mathbb{S}^{2} $.
\end{proposition}
Let $ (M,g) $ be a complete Riemannian 3-manifold with one end, finite first Betti number and $ R_{g}\geq 0 $. Let $ u $ be the solution to (\ref{u3}). We denote  
$$ 
\Sigma_{t}=\{x\in M| u(x)=1-\frac{1}{t}\},\,\,\,\Omega_{t}=\{p\}\cup\{x\in M| u(x)< 1-\frac{1}{t} \}
$$
Set
$$
\mathcal{T}=\Big\{t\in (0, +\infty)\,\,\Big|\,\, 1-\frac{1}{t} \text{ is a regular value of } u \Big\}
$$ 
Let $ \beta=b_{1}(M)<\infty $ be the first Betti number of $ M $. The following estimate is needed because the level sets $ \Sigma_{t} $ may be not connected.
\begin{lemma}\cite{You Z}\label{connect}
Let M be a connected manifold and let $ U\Subset M $ be a connected smooth domain such that $ M\setminus U $ is connected. If
$$
\partial U=\Sigma_{1}\cup\cdot \cdot \cdot \cup \Sigma_{N}
$$
then
$$
N-1\leq b_{1}(M).
$$
\end{lemma}
\begin{lemma}
For each regular  $ t\in \mathcal{T} $, the set $ \Omega_{t} $ is connected and $M\setminus \Omega_{t}$ is connected.
Furthermore, if $ N(t) $ is the number of components of $ \Sigma_{t} $, then
\begin{equation}\label{Nt}
N(t)-1\leq b_{1}(M).
\end{equation}
\begin{proof}
First we claim that $ \Omega_{t} $ is connected for each $ t\in \mathcal{T} $.  In fact, if this is not true, then there exists a connected component $ \Omega_{t}^{\prime} $ of $ \Omega_{t} $ which does not contain the pole $ p $ of $ u $. Moreover, $ \Omega_{t}^{\prime} $ is bounded by the fact
that $ \lim\limits_{x\to +\infty}u(x)=1 $. Then the harmonic function $ u $ on $ \Omega_{t}^{\prime} $ must achieve its minimum
in the interior of $ \Omega_{t}^{\prime} $, which is a contradiction.

Next, we claim that $M\setminus \Omega_{t}$ contains no bounded components for each $ t\in \mathcal{T} $. Otherwise, on such a bounded component the function $ u $ would achieve its maximum
at an interior point, which is again a contradiction. Since $ M $ has only one end, $M\setminus \Omega_{t}$ has only one unique unbounded component. Hence $M\setminus \Omega_{t}$ is connected for every $ t\in \mathcal{T} $. Now apply Lemma~\ref{connect} to obtain (\ref{Nt}).
\end{proof}
\end{lemma}

Along the level set $ \Sigma_{t} $, we define 
$$
\mathcal{A}(t)=t^{2}\int_{\Sigma_{t}}|\nabla u|^{2},
$$
$$
\mathcal{G}(t)=-\frac{4\pi}{t}+t\int_{\Sigma_{t}}|\nabla u|^{2},
$$
$$
\mathcal{F}(t)=4\pi t+t^{3}\int_{\Sigma_{t}}|\nabla u|^{2}-t^{2}\int_{\Sigma_{t}}|\nabla u|H.
$$
As in \cite{Miao}, we introduce the $ (0,2) $-symmetric tensors of $ M $:
$$
T=\nabla^{2}u+\frac{3u}{1-u}du\otimes du,\,\,\,B=T-\frac{tr(T)}{3}g,
$$
We define
$$
B_{1}(t)=\int_{\Sigma_{t}}\frac{|B|^{2}}{|\nabla u|^{2}},\,\,\,R_{1}(t)=\int_{\Sigma_{t}}R_{g}.
$$
Observe that 
$$
\frac{|B|^{2}}{2|\nabla u|^{2}}=\frac{|\nabla^{\Sigma_{t}}|\nabla u| |^{2}}{|\nabla u|^{2}}+\frac{1}{2}|\overset{\circ}{\mathbb{II}}|^{2}+\frac{3}{4}\Big(H-\frac{2}{1-u}|\nabla u|\Big)^{2}.
$$
By the work in Agostiniani-Mazzieri-Oronzio \cite{Ag3}, at a regular value $ t\in \mathcal{T} $, we have 
\begin{equation}\label{Fk2}
\begin{aligned}
\mathcal{F}^{\prime}(t)&=4\pi-2\pi\chi(\Sigma_{t})+\frac{1}{2}\int_{\Sigma_{t}}\Big(\frac{|B|^{2}}{|\nabla u|^{2}}+R_{g} \Big).
\end{aligned}
\end{equation}
By \cite{Wang1}, we know that there exists a $ \bar{t}>0 $ such that the level sets $ \Sigma_{t} $ are connected for all $ t\geq \bar{t} $. Hence, if the Proposition~\ref{YZ} does not hold, then there exists a regular $ t_{*}>\bar{t} $ such that, for every regular value $ t\geq t_{*} $, the level set $ \Sigma_{t} $ is a smooth, closed and connected surface of genus at least one. Therefore, by (\ref{Nt}), one has 
\begin{equation}\label{F 0}
\mathcal{F}^{\prime}(t)\geq 4\pi-4\pi(\beta+1)+ \frac{1}{2}(B_{1}(t)+R_{1}(t))\geq -4\pi \beta
\end{equation}
for all regular $ t> 0 $, and 
\begin{equation}\label{F 4pi}
\mathcal{F}^{\prime}(t)\geq 4\pi+\frac{1}{2}(B_{1}(t)+R_{1}(t))\geq 4\pi
\end{equation}
for all regular $ t\geq t_{*} $.
Now, we set 
$$ 
A_{1}=4\pi(\beta+1) t_{*} 
$$ 
and 
$$
\widetilde{\mathcal{A}}(t)=\mathcal{A}(t)-\frac{1}{2t}A_{1}.
$$
As in Proposition~\ref{diff ineq}, one has the following differential inequalities.
\begin{lemma}\label{diff ineq2}
Let $ (M,g) $ be a complete Riemannian 3-manifold with one end, finite first Betti number and $ R_{g}\geq 0 $. Let $ u $ be the solution of (\ref{u3}). Then we have the following differential inequalities
\\
$ (i) $ For all $ t\in (t_{*}, +\infty) $,
$$
t\mathcal{A}^{\prime}(t)\geq \mathcal{A}(t)-\frac{A_{1}}{t}+\frac{1}{2t}\int_{0}^{t}(B_{1}(s)+R_{1}(s))ds
$$
holds everywhere.
\\
$ (ii) $ For all $ t\in (t_{*}, +\infty)$,
$$
t\widetilde{\mathcal{A}}^{\prime}(t)\geq \widetilde{\mathcal{A}}(t)+\frac{1}{2t}\int_{0}^{t}(B_{1}(s)+R_{1}(s))ds
$$
holds everywhere.
\end{lemma}
\begin{proof}
At a regular value $ t\in \mathcal{T} $, we compute
$$
\mathcal{A}^{\prime}(t)=2t\int_{\Sigma}|\nabla u|^{2}-\int_{\Sigma_{t}}|\nabla u|H.
$$
Then we get 
\begin{equation}\label{t Ak2}
\begin{aligned}
t\mathcal{A}^{\prime}(t)-\mathcal{A}(t)=&t^{2}\int_{\Sigma}|\nabla u|^{2}-t\int_{\Sigma_{t}}|\nabla u|H\\
=&\frac{1}{t}\mathcal{F}(t)-4\pi\\
=&\frac{1}{t}\int_{0}^{t} \mathcal{F}^{\prime}(s)ds-4\pi,
\end{aligned}
\end{equation}
where we have used the fact that $ \lim\limits_{t\to 0}\mathcal{F}(t)=0 $ in \cite{Ag3}. On the other hand, it follows from (\ref{F 0}) and (\ref{F 4pi}) that 
$$
\begin{aligned}
t\mathcal{A}^{\prime}(t)-\mathcal{A}(t)&\geq \frac{4\pi(t-t_{*})}{t}-\frac{4\pi\beta t^{*}}{t}-4\pi+\frac{1}{2t}\int_{0}^{t}(B_{1}(s)+R_{1}(s))ds\\
&=-\frac{4\pi(\beta+1) t_{*}}{t}+\frac{1}{2t}\int_{0}^{t}(B_{1}(s)+R_{1}(s))ds\\
&=-\frac{A_{1}}{t}+\frac{1}{2t}\int_{0}^{t}(B_{1}(s)+R_{1}(s))ds.
\end{aligned}
$$
By Sard's theorem and continuity of $ \mathcal{A}(t) $, the above inequality holds for every $ t\in (t_{*}, +\infty) $. This proves $ (i) $.

To prove $ (ii) $, a simple computation shows that
$$
\begin{aligned}
t\widetilde{\mathcal{A}}^{\prime}(t)-\widetilde{\mathcal{A}}(t)&=t\mathcal{A}^{\prime}(t)-\mathcal{A}(t)+\frac{1}{t}A_{1}\geq \frac{1}{2t}\int_{0}^{t}(B_{1}(s)+R_{1}(s))ds.
\end{aligned}
$$
This proves $ (ii) $.
\end{proof}
Arguing as in Proposition~\ref{quadratic growth}, we have the following result.
\begin{lemma}
	Under the assumptions in Lemma~\ref{diff ineq2}, if $ \mathcal{A}(t_{0})>\frac{A_{1}}{t_{0}} $ for some $ t_{0}\in (t_{*},+\infty ) $, then there exists a positive constant $ C>0 $ such that
$$
\mathcal{A}(t)\geq Ct^{2}.
$$
for every $ t\geq t_{0} $.
\end{lemma}
\textbf{Proof of Proposition 4.1}. Suppose that the Proposition~\ref{YZ} does not hold, then there exists a regular $ t_{*}>0 $ such that, for every regular value $ t\geq t_{*} $, the level set $ \Sigma_{t} $ is a smooth, closed and connected surface of genus at least one.  We claim that 
\begin{equation}\label{Ak A2}
\mathcal{A}(t)\leq\frac{A_{1}}{t}
\end{equation}
for all $ t\in (t_{*}, +\infty) $. In fact, if there exists a $ t_{0}\in (t_{*}, +\infty)  $ such that $ \mathcal{A}(t_{0})>\frac{A_{1}}{t_{0}} $, then 
$$ 
\mathcal{A}(t)\geq Ct^{2} 
$$ 
for all $ t\in [t_{0}, +\infty) $, , which can be rewritten as 
$$
\int_{\Sigma_{t}}|\nabla u|^{2}\geq C.
$$
However, by the Cheng-Yau gradient estimate \cite{Cheng-Yau}, 
$$
\sup_{\Sigma_{t}}|\nabla u|=\sup_{\Sigma_{t}}|\nabla (1-u)|\leq \frac{C}{t}\to 0
$$
as $ t\to +\infty $. Therefore, we get 
$$
\int_{\Sigma_{t}}|\nabla u|^{2}\leq \sup_{\Sigma_{t}}|\nabla u|\int_{\Sigma_{t}}|\nabla u|=\mathfrak{c}_{\Sigma}\sup_{\Sigma_{t}}|\nabla u|\to 0.
$$
This is a contradiction.  Hence (\ref{Ak A2}) holds, which implies that 
$$
\lim\limits_{t\to +\infty}t\int_{\Sigma_{t}}|\nabla u|^{2}=0.
$$
Hence, we obtain that 
\begin{equation}\label{G lim2}
\lim\limits_{t\to +\infty}\mathcal{G}(t)=\lim\limits_{t\to +\infty}\Big[-\frac{4\pi}{t}+t\int_{\Sigma_{t}}|\nabla u|^{2}\Big]=0.
\end{equation}
However, by the Corollary~1.3 in \cite{Ag3}, we have 
\begin{equation}\label{G 00}
\lim\limits_{t\to 0}\mathcal{G}(t)=0.
\end{equation}
On the other hand, at a regular value $ t\in \mathcal{T} $, we have 
\begin{equation}\label{G F}
\mathcal{G}^{\prime}(t)=\frac{1}{t^{3}}\mathcal{F}(t)\geq \lim\limits_{t\to 0}\mathcal{F}(t)=0.
\end{equation}
Thus $ \mathcal{G}(t) $ is monotone nondecreasing in $ t $. It follows from (\ref{G lim2}) and (\ref{G 00}) that 
$$
\mathcal{G}(t)\equiv0
$$
for all $ t\in (0,\infty) $. Hence by (\ref{G F}), we conclude that 
$ \mathcal{F}(t)\equiv0 $ and so that 
$$
\mathcal{F}^{\prime}(t)\equiv0
$$
for all $ t\in (0,\infty) $, which contradicts to (\ref{F 4pi}). This completes the proof of Proposition~\ref{YZ}.
$\hfill\square$

\end{document}